\documentclass[12pt]{amsart}
\usepackage{amssymb}
\usepackage{amsmath}
\usepackage{pdfsync}

\begin{document}
\numberwithin{equation}{section}

\setcounter{tocdepth}{1}

\newcommand{\tslash}{{\backslash\!\!\!t}}
\newcommand{\dslash}{{\backslash\!\!\!d}}
\newcommand{\parslash}{{\backslash\!\!\!\partial}}
\newcommand{\Dslash}{{\backslash\!\!\!\!D}}

\def\1#1{\overline{#1}}
\def\2#1{\widetilde{#1}}
\def\3#1{\widehat{#1}}
\def\4#1{\mathbb{#1}}
\def\5#1{\frak{#1}}
\def\6#1{{\mathcal{#1}}}
\def\7#1{{\bf{#1}}}

\def\C{{\4C}}
\def\R{{\4R}}
\def\N{{\4N}}
\def\Z{{\4Z}}

\title[]{A Compact Expression\\ for Periodic Instantons}
\author[S.A. Cherkis]{Sergey A. Cherkis*}
\author[C. O'Hara]{Clare O'Hara}
\author[D. Zaitsev]{Dmitri Zaitsev}
\address{Sergey A. Cherkis: \newline
Department of Mathematics, the University of Arizona, Tucson AZ 85721, USA}
\email{\tt cherkis@math.arizona.edu}
\address{ Clare O'Hara and Dmitri Zaitsev:\newline School of Mathematics, Trinity College Dublin, Dublin 2, Ireland}
\email{\tt clohara@tcd.ie, zaitsev@maths.tcd.ie}
\thanks{* Supported by a grant from the Simons Foundation \#245643.}

\begin{abstract}{Instantons on various spaces  can be constructed via a generalization of the Fourier transform called the ADHM-Nahm transform.  An explicit use of this construction, however, involves rather tedious calculations.  Here we derive a simple formula for instantons on a space with one periodic direction.  It simplifies the ADHM-Nahm machinery and can be generalized to other spaces.}
\end{abstract}
\maketitle


\def\Label#1{\label{#1}}


\def\cn{{\C^n}}
\def\cnn{{\C^{n'}}}
\def\ocn{\2{\C^n}}
\def\ocnn{\2{\C^{n'}}}


\def\dist{{\rm dist}}
\def\const{{\rm const}}
\def\rk{{\rm rank\,}}
\def\id{{\sf id}}
\def\tr{{\sf tr}}
\def\aut{{\sf aut}}
\def\Aut{{\sf Aut}}
\def\End{{\sf End}}
\def\CR{{\rm CR}}
\def\GL{{\sf GL}}
\def\Re{{\sf Re}\,}
\def\Im{{\sf Im}\,}
\def\span{\text{\rm span}}
\def\gl{{\5g\5l}}
\def\sll{{\5s\5l}}
\def\u{{\5u}}
\def\su{{\5s\5u}}
\def\supp{{\sf supp}}
\def\Vec{{\sf Vec\,}}

\def\codim{{\rm codim}}
\def\crd{\dim_{{\rm CR}}}
\def\crc{{\rm codim_{CR}}}

\def\phi{\varphi}
\def\eps{\varepsilon}
\def\d{\partial}
\def\a{\alpha}
\def\b{\beta}
\def\g{\gamma}
\def\G{\Gamma}
\def\D{\Delta}
\def\Om{\Omega}
\def\k{\kappa}
\def\l{\lambda}
\def\L{\Lambda}
\def\z{{\bar z}}
\def\w{{\bar w}}
\def\Z{{\1Z}}
\def\t{\tau}
\def\th{\theta}

\emergencystretch15pt
\frenchspacing

\newtheorem{Thm}{Theorem}[section]
\newtheorem{Cor}[Thm]{Corollary}
\newtheorem{Pro}[Thm]{Proposition}
\newtheorem{Lem}[Thm]{Lemma}

\theoremstyle{definition}\newtheorem{Def}[Thm]{Definition}

\theoremstyle{remark}
\newtheorem{Rem}[Thm]{Remark}
\newtheorem{Exa}[Thm]{Example}
\newtheorem{Exs}[Thm]{Examples}

\def\bl{\begin{Lem}}
\def\el{\end{Lem}}
\def\bp{\begin{Pro}}
\def\ep{\end{Pro}}
\def\bt{\begin{Thm}}
\def\et{\end{Thm}}
\def\bc{\begin{Cor}}
\def\ec{\end{Cor}}
\def\bd{\begin{Def}}
\def\ed{\end{Def}}
\def\br{\begin{Rem}}
\def\er{\end{Rem}}
\def\be{\begin{Exa}}
\def\ee{\end{Exa}}
\def\bpf{\begin{proof}}
\def\epf{\end{proof}}
\def\ben{\begin{enumerate}}
\def\een{\end{enumerate}}
\def\beq{\begin{equation}}
\def\eeq{\end{equation}}

\section{Introduction}
A connection with self-dual curvature on a hermitian vector bundle over $\mathbb{R}^3\times S^1$ with finite action is called a {\em periodic instanton} or a {\em caloron}.  The integral of the second Chern class of such a connection is called its {\em caloron number} or the {\em number of calorons}.  Such connections play a crucial role in understanding gauge theories at finite temperature \cite{Gross:1980br}.  The first such connections with trivial asymptotic  holonomy around the $S^1$ were constructed by Harrington and Shepard in \cite{Harrington:1978ve}.  A general construction for all calorons was discovered by Nahm in \cite{NahmCalorons}. However, applying the Nahm transform in order to obtain explicit solutions can be difficult.  It was applied to obtain the caloron of caloron number one for the gauge group $SU(2)$ with a generic holonomy at infinity in \cite{Kraan:1998pm,Kraan:1998kp} and \cite{Lee:1998bb}.  The results of Kraan and van Baal \cite{Kraan:1998pm,Kraan:1998kp} were generalized by them to the case of one $SU(n)$ caloron in \cite{Kraan:1998sn} and \cite{Kraan:1998gh}.  Here we extend these results to the case of arbitrary caloron number and gauge group $SU(n).$ 

The Nahm transform for $k$ calorons with the gauge group $SU(n)$  \cite{NahmCalorons} is a generalization of the Atiyah, Drinfeld, Hitchin, and Manin (ADHM) construction.  It amounts to 
\begin{itemize}
\item finding some $k\times k$ Nahm data (see \S\ref{data}) on a circle with $n$ marked points $\lambda_\alpha, \alpha=1,\ldots,n$,  satisfying certain nonlinear ordinary differential equations (Nahm's equations), 
\item finding the kernel of a family of Dirac operators $\7D_t^\dag$ constructed out of the Nahm data, and  
\item finding the induced connection on the kernel bundle.
\end{itemize}
In other words to obtain a caloron solution one has to 1) solve a nonlinear system of the Nahm equations
\begin{equation}
i\frac{d}{ds} T_j+[T_0, T_j]+[T_{j+1}, T_{j+2}]=0
\end{equation} 
with some matching conditions at $s=\lambda_\alpha,$
2) use the solution to form a rank $2k$ linear system of ordinary differential equations  $\7D^\dagger\Psi=0$ and solve it, and
3) if the columns of $\Psi$ form an orthonormal basis of solutions, compute integrals of the form $A_\mu=\int \Psi^\dagger\partial_\mu\Psi ds,$ obtaining the caloron connection $A=\sum_\mu A_\mu dx^\mu.$ 

The relations we derive here simplify this standard Nahm transform procedure by a half.  In particular, instead of 2) one has to find the Green's function $F(s,t)=(\7D^\dagger\7D)^{-1},$ which amounts to solving a rank $k$ system of ODEs, and we eliminate step 3) by writing a simple expression in terms of the values of the Green's function at the marked points $F(\lambda_\alpha, \lambda_{\beta}).$

We introduce our notation in the next section.   Section~\ref{oper} contains basic operator relations and invertability arguments. 
In Section~\ref{Solution} we derive our main result Eq.~\eqref{A'''} expressing the $SU(n)$ caloron, of arbitrary instanton number and no monopole charges, in terms of the Nahm Green's function values at $n$ points.

\section{Setup}\Label{Setup}
\subsection{Spinors}
By a {\em spinor bundle} over a circle $\4S^1$ we mean any hermitian complex vector bundle $S\to M$ of rank $2$ over a manifold $M$, together with an action of the quaternion algebra $\4H$ by antihermitian bundle endomorphisms. That is, for every $x\in M$, we have a representation of $\4H$
in the Lie algebra $\u(S_x)$ of all antihermitian endomorphisms.
Equivalently, we have antihermitian bundle endomorphisms $e_1,e_2,e_3$
corresponding to the quaternion units and satisfying the quaternion relations.
We shall adopt the notation that whenever a symbol $u_j$ is defined for $j=1,2,3$,
we extend it by setting $u_{j+3i}=u_j$ for any $i\ge1$.
Then the relations for $e_j$ can be written as
\begin{equation}\Label{quat-rel}
e_j^2=-\id, \quad 
e_j e_{j+1}= -e_{j+1} e_j = e_{j+2}, \quad j=1,2,3.
\end{equation}
We also set $e_0:=\id$ for the identity operator corresponding to the real unit $1\in \4H$.
Thus $e_0$ is hermitian  whereas $e_1,e_2,e_3$ are antihermitian
bundle endomorphisms of $S$. It follows that on each fibre of $S$, $e_\mu$ are automatically linearly independent and hence form a basis of the space of all endomorphisms of the fibre.
Alternatively we can use the imaginary multiples
$\sigma_{\mu}:=ie_{\mu}$, $\mu=0,1,2,3$.
Then $\sigma_0=i\id$ is antihermitian and $\sigma_1,\sigma_2,\sigma_3$ are hermitian 
and satisfy the relations
\begin{equation}\Label{sigma-rel}
\sigma_j^2=\id, \quad 
\sigma_j \sigma_{j+1}= -\sigma_{j+1} \sigma_j = i\sigma_{j+2}, \quad j=1,2,3.
\end{equation}
One choice of $\sigma_j$ is Pauli matrices:
\begin{equation}\Label{pauli}
\sigma_1= 
\begin{pmatrix}
0&1\cr1&0
\end{pmatrix}, \quad
\sigma_2= 
\begin{pmatrix}
0&-i\cr i&0
\end{pmatrix}, \quad
\sigma_3= 
\begin{pmatrix}
1&0\cr0&-1
\end{pmatrix}.
\end{equation}

In the following we adopt the convention that the Greek index $\mu$ runs through $0,1,2,3$,
while the Latin index $j$ runs through $1,2,3$.

Any complex bundle $S$ over the unit circle $\4S^1$ is topologically trivial, given a trivialization 
$S\cong \4S^1\times\C^2$, where $\C^2$ is identified with the space of quaternions $\4H$ 
such that the complex structure on $\C^2$ corresponds to the right quaternion multiplication by $i\in\4H$
and the action of $e_1,e_2,e_3$ corresponds to the left quaternion multiplication
by the units. Thus we shall consider $S= \4S^1\times\C^2$ as trivial bundle over $\4S^1$ with trivial connection.

\subsection{Nahm data}\Label{data}
Let $E\to \4S^1$ be a hermitian complex vector bundle of rank $k$ over the unit circle $\4S^1$.
We write $\langle\cdot|\cdot\rangle$ for the hermitian metric which is assumed
to be complex-antilinear in the first and complex-linear in the second argument.
We shall identify the points of $\4S^1$ with points in $\R$ via the exponential map.
That is, we simply write $s$ for $e^{is}\in \4S^1$.
The {\em Nahm data} on $E$ consists of a hermitian 
(i.e.\ preserving the metric) connection $\nabla$ on $E$,
three hermitian bundle endomorphisms $T_1,T_2,T_3$ of $E$,
a partition $\L=\{\lambda_1,\ldots,\lambda_n\}$, $0\le\lambda_1<\ldots<\lambda_n<2\pi$
of $[0,2\pi)$ (which is identified with $\4S^1$), 
and a boundary data given by (complex-)linear maps 
\begin{equation}\Label{bdry}
Q_\alpha^\dag\colon (S\otimes E)_{\l_\alpha}\to W_\alpha,\quad \alpha=1,\ldots,n,
\end{equation}
where $(S\otimes E)_{\l_\alpha}$ is the fibre of $S\otimes E$ at $s=\lambda_\alpha$ and $W_\alpha$ is an auxiliary complex vector space
with some fixed hermitian metric for each $\alpha$.

Given the partition $\L$, we introduce the following function classes on $\4S^1$.
Denote by $\6C_\L$ the class of all functions on $\4S^1$ that
are continuous on $\4S^1$ and smooth up to the
boundary on each open interval $(\lambda_{\alpha},\lambda_{\alpha+1})$, $\alpha=1,\ldots,n$,
where we adopt the notation 
\begin{equation}\Label{l-notation}
\lambda_{\alpha+nl}=\lambda_\alpha+2\pi l, \quad l=0,\pm1,\pm2,\ldots.
\end{equation}
By $\6C'_\L$ we denote the class of all functions on the union 
$\cup_{\alpha=1}^n (\lambda_{\alpha},\lambda_{\alpha+1})\subset\4S^1$
which are smooth up to the boundary on each $(\lambda_{\alpha},\lambda_{\alpha+1})$.
Thus a function from $\6C'_\L$ may have different limits at $s=\l_\alpha$
from the left and from the right and is not defined at $s=\l_\alpha$.
Finally $\6C''_\L$ denotes the class of functions of the form
$f(s)+\sum_{\alpha=1}^n \delta(s-\l_\alpha)a_\alpha$, where $f$ is any function in $\6C'_\L$,
$a_1,\ldots,a_\alpha$ are any coefficients and $\delta(s)$ is Dirac's delta-function.
Similarly, for any hermitian vector bundle $V$ on $\4S^1$, we define
spaces of sections of classes $\6C_\L$, $\6C'_\L$, $\6C''_\L$
that we denote respectively by $\G_\L(V)$, $\G'_\L(V)$, $\G''_\L(V)$.
The definitions of $\G_\L(V)$ and $\G'_\L(V)$ are straightforward, whereas
for $\G''_\L(V)$ we need to consider generalized sections $\delta(s-\l_\alpha)v$ with $v$
being any vector in the fiber $V_{\l_\alpha}$ of $V$.
We regard $\delta(x-\l_\alpha)v$ as the (complex-)linear functional on $\G_\L(V)$
assigning to every section $g\in \G_\L(V)$ the complex number $\langle v| g(\l_\alpha)\rangle$.
A section $f\in\G'_\L(V)$ also defines the linear functional given by the $L^2$ product
$$g\in\G_\L(V)\mapsto \int \langle f | g\rangle\, ds,$$
where the integral is always taken over $\4S^1$ with respect to the standard volume form $ds$ unless specified otherwise.
Thus the class of sections $\G''_\L(V)$ can be rigorously seen
as a subspace of the space of all linear functionals on $\G_\L(V)$.
As is customary, we also write the functional induced by $\delta(x-\l_\alpha)v$ in the form
$$g\mapsto \langle v| g(\l_\alpha)\rangle  = \int \langle \delta(s-\l_\alpha)v|g(s)\rangle\, ds.$$

It is straightforward to define connections and bundle endomorphisms of classes $\6C_\L$ and $\6C'_\L$.
We assume that $\nabla$ and the $T_j$ are of class $\6C'_\L$.
In particular, this means that in a (local) trivialization $E|_{(a,b)}\cong (a,b)\times\C^k$ 
on an interval $(a,b)\subset\4S^1$ (with metric also trivialized), we can write 
\begin{equation}\Label{nabla-loc}
\nabla = \frac{d}{ds} -iT_0,
\end{equation}
where $s$ is the coordinate in the interval $(a,b)$ and $T_0$ is an hermitian $k \times k$ matrix valued function of class $\6C'_\L$.

\subsection{Weyl operator}
For the sake of simpler notation, we shall write $t$ instead of $t\id$ 
for the multiplication operator. 
For every $4$-vector $t=(t_0,t_1,t_2,t_3)\in\R^4$,
define the operator
\begin{equation}\Label{D}
\7D=\7D_t\colon \G_\L(S\otimes E)\to \G'_\L(S\otimes E),
\quad \7D:=\id\otimes (i\nabla +t_0) + \sum_j e_j\otimes (T_j+t_j),
\end{equation}
where the Latin index $j$ runs through $1,2,3$ according to our convention.
Using the induced hermitian metric $\langle \cdot|\cdot \rangle$ on $S\otimes E$,
we consider the {\em adjoint} operator $\7D^\dag$ given by
\begin{equation}\Label{adjoint}
\7D^\dag =\7D^\dag_t\colon \G'_\L(S\otimes E)\to \G''_\L(S\otimes E),
\int \langle \7D^\dag f| g\rangle\, ds = \int\langle f |\7Dg \rangle\, ds,
\end{equation}
with $ f\in \G'_\L(S\otimes E), g\in \G_\L(S\otimes E).$
That is, the result of applying $\7D^\dag$ to $f\in \G'_\L(S\otimes E)$ is the functional on 
$ \G_\L(S\otimes E)$ given by 
$$g\in\G_\L(S\otimes E) \mapsto  \int\langle f |\7Dg \rangle\, ds.$$
Note that since the section $f$ is allowed to have discontinuities on $\L$,
its derivative may not correspond to any `true' section (not even in $L^2$) 
but is, on the other hand, always realized as a functional on $\G_\L(S\otimes E)$
(involving delta-functions).

The fact that the connection $\nabla$ and operators $\sigma_j$ and $T_j$ are hermitian (and hence $e_j$ are anti-hermitian) implies (via integration by parts) that 
\begin{equation}\Label{D-adj}
\7D^\dag=\id\otimes(i\nabla+t_0) - \sum_j e_j\otimes (T_j+t_j),
\end{equation}
where the connection operator $\nabla$ sends sections of class $\6C'_\L$ to sections of class $\6C''_\L$
(seen as functionals) according to the formula
$$\int \langle \nabla f | g\rangle\, ds = -\int \langle f | \nabla g\rangle\, ds, \quad g\in \G_\L(E).$$
Note that it used here that the connection $\nabla$ is hermitian, i.e.\ 
\begin{equation}\Label{herm}
\langle \nabla f | g\rangle + \langle f|  \nabla g\rangle = d\langle f | g\rangle
\end{equation}
holds for smooth sections $f$ and $g$.

We next write $W_\L:=\oplus_{\alpha=1}^n W_\alpha$ for the boundary data \eqref{bdry} 
and consider the boundary data operator
\begin{equation}\Label{Qdag}
\7Q^\dag\colon \G_\L(S\otimes E) \to W_\L,
\quad \7Q^\dag f= \big( Q^\dag_1f(\l_1),\ldots,Q^\dag_nf(\l_n)).
\end{equation}
We equip $W_\L$ with the direct sum of hermitian metrics on $W_\alpha$
and consider the adjoint
$$\7Q\colon W_\L\to \G''_\L(S\otimes E)$$
given by
$$\quad \int \langle \7Q v |g\rangle\, ds = \langle v | \7Q^\dag g\rangle, \quad g\in\G_\L(S\otimes E).$$
Explicitly, we have
\begin{equation}\Label{Q}
\7Q(v_1,\ldots,v_n)=\sum_{\alpha=1}^n \delta(s-\l_\alpha)Q_\alpha v_\alpha,
\end{equation}
where $Q_\alpha\colon W_\alpha\to (S\otimes E)_{\l_\alpha}$ is the adjoint of $Q_\alpha^\dag$.

The {\em Weyl operator} is now defined by
\begin{equation}\Label{weyl}
\6D :=  \7D \oplus \7Q^\dag\colon \G_\L(S\otimes E) \to \G'_\L(S\otimes E)\oplus W_\L
\end{equation}
and its adjoint is
$$\6D^\dag =  \7D^\dag + \7Q\colon \G'_\L(S\otimes E)\oplus W_\L \to \G''_\L(S\otimes E),$$
where $\7D^\dag$ acts on the first component of the direct sum and $\7Q$ acts on the second component.

\section{Operator operations}\label{oper}

\subsection{Composing operators with their adjoints}
Using the relations \eqref{sigma-rel}
we have
\begin{multline}\Label{dd}
\7D^\dag \7D =
\big(\id\otimes (i\nabla+t_0) - \sum_j e_j\otimes (T_j+t_j)\big)
\big(\id\otimes (i\nabla+t_0) + \sum_le_l\otimes (T_l+t_l)\big)\\
=
\id\otimes\big((i\nabla+t_0)^2 + \sum_j(T_j+t_j)^2\big)
 +\sum_j e_j \otimes \big([i\nabla, T_j] + [T_{j+1},T_{j+2}] \big),
\end{multline}
which is a second order differential operator from $\G_\L(S\otimes E)$ into $\G''_\L(S\otimes E)$. Note that $[i\nabla, T_j]$ is a zero order operator of class $\6C''_\L$
acting by a bundle endomorphism.
It is the content of the {\em Nahm Equations} that the operators
\beq
[i\nabla, T_j] + [T_{j+1},T_{j+2}]
\eeq
in the last term vanish on each interval $(\l_\alpha,\l_{\alpha+1})$. 
However, they still contribute with delta-functions supported on $\L$,
corresponding to the jumps of $T_j$.
The latter are being taken care of by the boundary operator $\7Q$.

In view of \eqref{Qdag} and \eqref{Q} we have
\begin{equation}
\7Q \7Q^\dag f = \sum_{\alpha=1}^{n}\delta(s-\l_{\alpha})Q_{\alpha}Q^{\dag}_{\alpha}f(\l_{\alpha}).
\end{equation}
Since $e_0,e_1,e_2,e_3$ is a basis over $\C$ 
in the space of all endomorphisms of $S$, we can write
\begin{equation}\Label{qq}
\7Q \7Q^\dag= \sum_\mu e_\mu\otimes(\7Q \7Q^\dag)_\mu,
\end{equation}
where the summation over Greek indices runs over $0,1,2,3$, as agreed
and $(\7Q \7Q^\dag)_\mu$ are the corresponding components
acting from $\G_\L(E)$ into $\G''_\L(E)$.
The component $(\7Q \7Q^\dag)_0$ can be obtained directly by taking
the trace $\tr_S$ with respect to $S.$
Here, for any complex vector space $V$ and an endomorphism $A$ of $S\otimes V\cong V\oplus V$,
we write 
\begin{equation}\Label{tr-vec}
A=
\begin{pmatrix}
A_{11} & A_{12}\cr
A_{21} & A_{22}
\end{pmatrix}, \quad
\tr_{S} A = A_{11}+A_{22}, \quad 
\Vec A = 
A - \frac12\id_S\otimes\tr_S A.
\end{equation}
In fact, it follows from the relations \eqref{quat-rel}
that the trace of $e_j$ is zero for all $j=1,2,3$.
Then taking $\tr_{S}$ of both sides in \eqref{qq}
and using the fact that $e_0=\id$, we obtain
\begin{equation}\Label{trace}
\tr_S \7Q \7Q^\dag = 2 (\7Q \7Q^\dag)_0.
\end{equation}
Putting \eqref{dd} and \eqref{qq} together and using \eqref{trace} we obtain
\begin{multline}\Label{ddqq}
\6D^\dag \6D = (\7D^\dag+\7Q)(\7D+\7Q^\dag)\\
=\id\otimes\big((i\nabla+t_0)^2 + \sum_j(T_j+t_j)^2 + \frac12\tr_S(\7Q \7Q^\dag)\big)\\
 +\sum_j e_j \otimes \big([i\nabla, T_j] + [T_{j+1},T_{j+2}] + (\7Q \7Q^\dag)_j\big).
\end{multline}
Now we impose the following conditions that complement the Nahm Equations with boundary equations:
\begin{equation}\Label{nahm}
[i\nabla, T_j] + [T_{j+1},T_{j+2}] + (\7Q \7Q^\dag)_j=0.
\end{equation}
Thus \eqref{ddqq} becomes 
\begin{equation}\Label{ddqq'}
\6D^\dag \6D
=\id\otimes\big((i\nabla+t_0)^2 + \sum_j(T_j+t_j)^2 + \frac12\tr(\7Q \7Q^\dag)\big)
\end{equation}
and \eqref{dd} can be written as
\begin{equation}\Label{dd'}
\7D^\dag \7D =
\id\otimes\big((i\nabla+t_0)^2 + \sum_j(T_j+t_j)^2\big)
 -\sum_j e_j\otimes(\7Q\7Q^\dagger)_j.
\end{equation}

\subsection{Positivity and kernel of $\6D^\dag \6D$}
Formula \eqref{ddqq'} expresses $\6D^\dag \6D$
as a sum of the non-negative operators 
$$\6D^\dag \6D=A_{0}^{\dag}A_{0} + \sum_{j=1}^{3}A_{j}^{\dag}A_{j} + \sum_{m=4}^{5} A^{\dag}_{m}A_{m},$$
where
$$A_{0}:= \id\otimes(i\nabla+t_0), \quad A_{j}=\id\otimes (T_{j}+t_{j}),$$
and $\7Q$ is written as $\7Q=(A^{\dag}_{4},A^{\dag}_{5})$
acting by left matrix multiplication on
a section $f\in \G_{\L}(S\otimes E)\cong \G_{\L}(E)\oplus \G_{\L}(E)$
written as colum $f_{1}\choose {f_{2}}$ with $f_{1},f_{2}\in \G_{\L}(E)$.
The non-negativity here means
\begin{multline}
\int\langle f | A^{\dag}_{m}A_{m}f\rangle \, ds  =  \int \langle A_{m}f | A_{m}f\rangle \, ds\ge 0,\\
f\in\G_{\L}(S\otimes E), \quad m=0,\ldots, 5.
\end{multline}
In particular, a section $f\in\G_{\L}(S\otimes E)$
is in the kernel of $\6D^\dag \6D$ if and only if
it is in the kernel of each $A_{m}$, $m=0,\ldots,5$.
In the latter case $f$ has to be parallel 
with respect to the connection $\id\otimes(i\nabla + t_{0})$.
Hence, if $f\ne0$, then $f(0)\ne0$ has to be 
an eigenvector with eigenvalue $1$ of the monodromy operator
$\imath\colon (S\otimes E)_{0}\to (S\otimes E)_{0}$
obtained by following sections parallel with respect to $\id\otimes(i\nabla + t_{0})$
around $\4S^{1}$ (see \S\ref{ker} below for more details).
If $\imath_{0}$ is the monodromy for $t_{0}=0$,
then $\imath=e^{-2\pi t_{0}}\imath_{0}$ (see also \S\ref{ker}).
Therefore $\imath$ can only have eigenvalue $1$ for a discrete set of $t_{0}$.

Furthermore, if $f\ne0$ is in the kernel of $A_{j}$, $j=1,2,3$,
the scalar $t_{j}$ has to be an eigenvalue of $T_{j}$ for each $j$
(at every point $s$ where $f$ is not zero).
The latter is clearly possible only for finitely many values of $t_{j}$.
Summarizing, we conclude that $\6D^\dag \6D$
and hence $\6D$ has zero kernel for all $t=(t_{0},t_{1},t_{2},t_{3})$
except $t_{0}$ in a discrete set and each $t_{j}$ in a finite set.

\subsection{The kernel of $\7D^\dag$}\Label{ker}
We next consider the kernel of $\7D^\dag$.
In view of the formula \eqref{D-adj}, $\7D^\dag$ is
a first order linear differential operator
whose coefficients are smooth 
up to the boundary on each interval $(\l_\alpha,\l_{\alpha+1})$.
Hence 
\begin{equation}\Label{kernel}
\7D^\dag f=0
\end{equation}
 is a first order linear ODE system.
 It is well-known that any 
 local solution $f$ of such a system on a subinterval of $(\l_\alpha,\l_{\alpha+1})$
extends uniquely to a global solution on $(\l_\alpha,\l_{\alpha+1})$
which is smooth up to the boundary.
The solutions  are parametrized 
 by the initial boundary value $f(\l_\alpha)$
 (or, equivalently, by the value at 
any fixed point $s\in [\l_\alpha,\l_{\alpha+1}]$).
This parametrization is given by the smooth family of maps
\begin{equation}\Label{param}
v_{s_0}\colon f|_{s=s_0}\mapsto f \in \G([\l_\alpha,\l_{\alpha+1}],S\otimes E)
\end{equation}
such that $v_{s_0}(v_0)$ is a solution of \eqref{kernel}
with value $v_0$ in the fibre $(S\otimes E)_{s_0}$ over $s_0$.
Every solution arises in this way.
The smoothness of the family means that $v_{s_0}(f_0)(s)$
is smooth in $(s_0,s,f_0)$.
For every pair $s_0,s_1\in [\l_\alpha,\l_{\alpha+1}]$,
\eqref{param} defines a linear isomorphism 
\begin{equation}\Label{mon}
\imath_{s_1,s_0}\colon (S\otimes E)_{s_0} \to (S\otimes E)_{s_1}, \quad 
\imath_{s_1,s_0}(f_0):=v_{s_0}(f_0)(s_1).
\end{equation}
In particular, $\imath_{\l_{\alpha+1},\l_\alpha}\colon (S\otimes E)_{\l_\alpha} \to (S\otimes E)_{\l_{\alpha+1}}$
is the {\em monodromy} of $\7D^\dag$ over the interval $[\l_\alpha,\l_{\alpha+1}]$.

We now wish to determine solutions of $\eqref{kernel}$ on the whole $\4S^1$.
Note that by our definition, $\7D^\dag$ acts on $\G'_\L(S\otimes E)$.
Furthermore, it is clear from the formula \eqref{D-adj} that
any solution of \eqref{kernel} must be continuous
and therefore of class $\G_\L(S\otimes E)$.
Thus global solutions of \eqref{kernel} are obtained as concatenations
of its solutions over the intervals $(\l_\alpha,\l_{\alpha+1})$.
That is, we extend the construction of \eqref{mon}
to any $s_0\in [\l_{\alpha_0},\l_{\alpha_0+1})$ and $s_1\in [\l_{\alpha_1},\l_{\alpha_1+1})$ with $\alpha_0\le \alpha_1$
by concatenation of 
the monodromies \eqref{mon} over all intervals beginning at $s_0$
and ending at $s_1$:
\begin{equation}\Label{glue}
\imath_{s_1,s_0}:=\imath_{s_1,\l_{\alpha_1}}\circ \imath_{\l_{\alpha_1},\l_{\alpha_1-1}} \circ \ldots  \circ
 \imath_{\l_{\alpha_0+2},\l_{\alpha_0+1}} \circ \imath_{\l_{\alpha_0+1},s_0}.
\end{equation}
Then any solution $f$ of \eqref{kernel} over $(s_0,s_1)$
is of the form $f(s)=\imath_{s,s_0}(v_0)$ for some initial value
$v_0=f(s_0)\in (S\otimes E)_{s_0}$ and vice versa,
this formula gives a solution over $(s_0,s_1)$ for each $v_0\in  (S\otimes E)_{s_0}$.
If we now take $s_1=s_0+2\pi$, we see that the solutions of \eqref{kernel}
over $\4S^1$ correspond precisely to the initial data $v_0\in  (S\otimes E)_{s_0}$
such that $\imath_{s_0+2\pi,s_0}(v_0)=v_0$, i.e.\
to the eigenvectors of the monodromy operator $\imath_{s_0+2\pi,s_0}$
with eigenvalue $1$.
Summarizing, we have constructed for every operator $\7D^\dag=\7D^\dag_t$
the family of its monodromy operators $\imath_{s_0+2\pi,s_0}$
such that the operator has a nontrivial kernel if and only if
the monodromy operator has eigenvalue $1$ for any fixed $s_0$.
In particular, $\imath_{s_0+2\pi,s_0}$ has eigenvalue $1$
if and only $\imath_{s_1+2\pi,s_1}$ does.
(In fact, monodromies for different $s_0$
are conjugate to each other and therefore have equal eigenvalues.)

We now compare the eigenvalues of monodromies of $\7D^\dag_t$
for different $t$. It follows from the formula \eqref{D-adj} that 
\begin{equation}\Label{dt}
\7D^\dag_t=\7D^\dag_0 + t, \quad t:= t_0 - \sum_j e_j\otimes t_j,
\end{equation}
where the operator $t$ has constant coefficients
(in any local trivialization of $E$).
Assume that $f(s)$ is any solution of $\7D^\dag_{(0,t_1,t_2,t_3)}f=0$
defined on an interval $I$. Then it follows from the formula
\begin{equation}\Label{t0}
\7D^\dag_t e^{it_0s}f(s) =e^{it_0s}(\7D^\dag_t - t_0)f(s))
\end{equation}
that $\2f(s):=e^{it_0s}f(s)$ is a solution of $\7D^\dag_t \2f=0$ on $I$.
Hence the monodromies $\imath_{s_0+2\pi,s_0}$ and $\2\imath_{s_0+2\pi,s_0}$
of $\7D^\dag_{(0,t_1,t_2,t_3)}$ and $\7D^\dag_t$ respectively are related by
\begin{equation}\Label{mon-relation}
\2\imath_{s_0+2\pi,s_0}= e^{2i\pi t_0}\imath_{s_0+2\pi,s_0}.
\end{equation}
In particular, the monodromy $\2\imath_{s_0+2\pi,s_0}$ has eigenvalue $1$
only for a discrete (but nonempty) set of values $t_0$.
Furthermore, $\7D^\dag_t$ depends analytically on $t$.
Then it follows from the construction of the monodromy
that the monodromy also depends analytically on $t$.
Therefore the monodromy has eigenvalue $1$ precisely
for all $t$ in a (real-)analytic subset $A\subset\R^4_{t_0,t_1,t_2,t_3}$.
The above conclusion implies that $A$ is a proper analytic subset of $\R^4$
and has discrete (and nonempty) intersection with each line $\R_{t_0}\times\{(t_1,t_2,t_3)\}$.
Summarizing, we conclude that $\7D^\dag$ has zero kernel for all $t\in\R^4$
away from a proper analytic subset.

A completely analogous construction can be done for 
$\7D\colon \G_\L(S\otimes E)\to \G'_\L(S\otimes E)$.
In particular, $\7D=\7D_t$ has zero kernel for all $t$
outside a proper analytic subset of $\R^4$.

\subsection{The Fundamental Solution for $\7D^\dag$}\Label{fundamental}
We now assume that $t$ is chosen such that $\7D^\dag$ has zero kernel.
In this case we construct the inverse of $\7D^\dag$ by finding its fundamental solution as follows.
Recall that a {\em fundamental solution} for $\7D^\dag$  is a family of linear operators
$$B(x,y)\colon (S\otimes E)_y\to (S\otimes E)_x, \quad x,y\in\4S^1,$$
safisfying
\begin{equation}\Label{green-rel}
\7D^\dag B(\cdot,y) = \delta(\cdot-y),
\end{equation}
where the $\delta$-function is seen as a linear operator
from $(S\otimes E)_y$ into $\G''_\L(S\otimes E)$,
defined in the obvious way: $v\mapsto \delta(\cdot-y)v$.
To construct the fundamental solution,  consider the monodromy operator
$$\imath_{s_1,s_0}\colon (S\otimes E)_{s_0}\to (S\otimes E)_{s_1}.$$
Then  
\begin{equation}\Label{parf}
f(s)=\imath_{s,s_0}v_0
\end{equation}
parametrizes the solutions of $\7D^\dag f=0$
for $v_0\in (S\otimes E)_{s_0}$ which we consider for $s\in [s_0,s_0+2\pi)$.
Then it follows from the formula \eqref{D-adj}
that $f$ satisfies
\begin{equation}\Label{del-eq}
\7D^\dag f (s) = \delta(s-s_0)v
\end{equation}
for some $v\in (S\otimes E)_{s_0}$
if and only if 
\begin{equation}\Label{v}
\imath_{s_0+2\pi,s_0}v_0-v_0 = v,
\end{equation}
i.e.\ $v$ is the discontinuity of $f$ at $s_0$.
Note that here $f\in\G_{\L\cup\{s_0\}}(S\otimes E)$
and therefore we have to extend $\7D^\dag$ to the operator
from $\G_{\L\cup\{s_0\}}(S\otimes E)$ into $\G''_{\L\cup\{s_0\}}(S\otimes E)$,
which is achieved by replacing $\L$ with $\L\cup\{s_0\}$ in the 
definition of $\7D^\dag$.

We now use our basic assumption that $\7D^\dag$ has zero kernel,
which is equivalent to $\imath_{s_0+2\pi,s_0}$ not having eigenvalue $1$.
Then the operator $\imath_{s_0+2\pi,s_0}-\id$ is invertible and hence
\eqref{v} can be solved in the form $v_0=(\imath_{s_0+2\pi,s_0}-\id)^{-1}v$.
Substituting into \eqref{parf} we obtain the section
\begin{equation}
f(s)=\imath_{s,s_0}(\imath_{s_0+2\pi,s_0}-\id)^{-1}v
\end{equation}
that solves \eqref{del-eq}. Setting
\begin{equation}
B(x,y):=\imath_{x,y}\circ(\imath_{y+2\pi,y}-\id)^{-1}
\end{equation}
we obtain the basic relation \eqref{green-rel} 
that shows that $B(x,y)$ is a fundamental solution (Green's function) for $\7D^\dag$.
It follows from the construction of the monondromy that
$B(x,y)$ is defined and continuous away from the diagonal $x=y$ and
is smooth away from
the lines $x=\l_\alpha$ and $y=\l_\alpha$.
The diagonal and these lines cut $\4S^1\times\4S^1$
into a union of triangles and rectangles and the fundamental solution  $B(x,y)$ extends smoothly
to the closure of each triangle. The main property of the  fundamental solution  
is that the corresponding integral operator
\begin{equation}\Label{B}
(\7B f)(x):=\int B(x,y) f(y)\, dy
\end{equation}
gives a solution of $\7D^\dag g=f$ in the form $g=\7Bf$.
Since $\7D^\dag$ maps $\G'_\L(S\otimes E)$ into $\G''_\L(S\otimes E)$,
we consider here $f\in \G''_\L(S\otimes E)$.
Then the integral in \eqref{B} is defined for $x$ in each interval $(\l_\alpha,\l_{\alpha+1})$
and is smooth in $x$ up to the boundary of this interval.
Thus $\7B f$ is a section of class $\6C'_\L$.
Summarizing we conclude that $\7B$ maps $\G''_\L(S\otimes E)$ into $\G'_\L(S\otimes E)$
and is a right inverse of $\7D^\dag$, i.e.\ $\7D^\dag\circ \7B = \id$.
In particular, $\7D^\dag$ is surjective onto $\G''_\L(S\otimes E)$.
Since $\7D^\dag$ has zero kernel by our assumption, it is also injective.
Therefore, $\7B$ is actually the (two-sided) inverse of $\7D^\dag$, thus we write $\7B=\big(\7D^\dagger\big)^{-1}.$

\subsection{Green's functions for $\7D^\dag\7D$ and $\6D^\dag\6D$}
Here we extend the above construction to the second order
differential operators $\7D^\dag\7D$ and $\6D^\dag\6D$
given by the formulae \eqref{dd'} and \eqref{ddqq'} respectively.
We shall assume that both 
$\7D=\7D_t$ and $\7D^\dag=\7D^\dag_t$ have zero kernel,
which holds for all $t$ outside a proper analytic subset of $\R^4$
in view of \S\ref{ker}.
Consequently,  the composition operator $\7D^\dag\7D$ acting from $\G_\L(S\otimes E)$
into $\G''_\L(S\otimes E)$ also has zero kernel.

We next construct the monodromy operator for $\7D^\dag\7D$.
Since it is a linear second order differential operator with smooth coefficients
on each $(\l_\alpha,\l_{\alpha+1})$, the solutions of the equation
\begin{equation}\Label{dd-ker}
\7D^\dag\7D f =0
\end{equation}
 on each inverval
are parametrized by their values $f(s)\in (S\otimes E)_s$ and the first order derivatives 
$(\nabla f)(s)\in (S\otimes E)_s$ at a fixed point $s\in [\l_\alpha,\l_{\alpha+1}]$.
Following the solutions of \eqref{dd-ker} as before, we define for $s_0,s_1\in [\l_\alpha,\l_{\alpha+1}]$,
the monodromy operator
\begin{equation}
\imath_{s_1,s_0}\colon (S\otimes E)_{s_0}\oplus (S\otimes E)_{s_0}
 \to (S\otimes E)_{s_1}\oplus (S\otimes E)_{s_1},
\end{equation}
where $f(s)\oplus g(s)=\imath_{s,s_0}(v_0,v'_0)$
if and only if $f(s)$ is the solution of \eqref{dd-ker} 
with initial data $v_0=f(s_0)$, $v'_0=(\nabla f)(s_0)$,
and $g(s)=(\nabla f)(s)$.

As the next step we glue the monodromies like in the formula \eqref{glue} above.
There is a new ingredient, however, due to the fact that the coefficients
of $\7D^\dag \7D$ may involve $\delta$-functions supported on $\L$.
Indeed, the latter may arise from the last term in \eqref{dd'}
as well as from $(i\nabla+t_0)^2$ which, in a local trivialization of $E$,
has in view of \eqref{nabla-loc} the form
\begin{equation}
\Big(i\frac{d}{ds}+T_0+t_0\Big)^2 = -\frac{d^2}{ds^2} 
+i\frac{dT_0}{ds}+(T_0+t_0)^2.
\end{equation}
Hence, in order for $f$ to be a solution of \eqref{dd-ker}
across the jumping point $\l_\alpha$, it has to be continuous at $\l_\alpha$
and the jump of its derivative at $\l_\alpha$ has to be equal to a certain linear
function of its value $f(\l_\alpha)$, which is determined by the $\delta$-function
terms of $\7D^\dag \7D$. Thus, for each $\l_\alpha$, we obtain
an additional monodromy operator 
\begin{equation}
\imath_{\l_\alpha}\colon (S\otimes E)_{\l_\alpha}\oplus (S\otimes E)_{\l_\alpha}
 \to (S\otimes E)_{\l_\alpha}\oplus (S\otimes E)_{\l_\alpha},
\end{equation}
such that $f$ is a solution of \eqref{dd-ker} in an interval $(s_0,s_1)$
with $\l_{\alpha-1}\le s_0<\l_\alpha<s_1\le \l_{\alpha+1}$ if and only if, for some $g$, $v_0$, $v'_0$,
$$(f(s),g(s))= 
\begin{cases}
\imath_{s,s_0}(v_0,v'_0) & s<\l_\alpha \\
\imath_{s,\l_\alpha}\circ \imath_{\l_\alpha} \circ \imath_{\l_\alpha,s_0} & s>\l_\alpha.
\end{cases}$$
Then the gluing formula for general $\imath_{s_1,s_0}$ with 
$s_0\in [\l_{\alpha_0},\l_{\alpha_0+1})$, $s_1\in [\l_{\alpha_1},\l_{\alpha_1+1})$ is given by
\begin{equation}\Label{glue'}
\imath_{s_1,s_0}:=\imath_{s_1,\l_{\alpha_1}}\circ \imath_{\l_{\alpha_1}} \circ \imath_{\l_{\alpha_1},\l_{\alpha_1-1}}\circ \imath_{\l_{\alpha_1-1}} \circ \ldots \circ \imath_{\l_{\alpha_0+2}} \circ
 \imath_{\l_{\alpha_0+2},\l_{\alpha_0+1}}\circ \imath_{\l_{\alpha_0+1}} \circ \imath_{\l_{\alpha_0+1},s_0},
\end{equation}
and for any $s_0\in\4S^1$, the solutions of \eqref{dd-ker} on the whole $\4S^1$
correspond precisely to the initial data 
$(v_0,v'_0)\in  (S\otimes E)_{s_0}\oplus (S\otimes E)_{s_0}$
such that $\imath_{s_0+2\pi,s_0}(v_0,v'_0)=(v_0,v'_0)$,
i.e.\ $(v_0,v'_0)$ is an eigenvector of $\imath_{s_0+2\pi,s_0}$
with eigenvalue $1$. Since the kernel of $\7D^\dag\7D$ is zero by our assumption,
the operator $\imath_{s_0+2\pi,s_0}-\id$ is invertible.
Furthermore, solutions of
\begin{equation}
\7D^\dag\7D f(s)=\delta(s-s_0)v, \quad v\in (S\otimes E)_{s_0}
\end{equation}
correspond to the initial data $(v_0,v'_0)$ satisfying
\begin{equation}
\imath_{s_0+2\pi,s_0}(v_0,v'_0)=(v_0,v'_0+v),
\end{equation}
or, equivalently, $(v_0,v'_0)=(\imath_{s_0+2\pi,s_0}-\id)^{-1}(0,v)$.
We can now write the Green's function, 
\begin{equation}
G(x,y):=\pi_2\circ\imath_{x,y}\circ (\imath_{y+2\pi,y}-\id)^{-1}\circ\pi_1,
\end{equation}
where we have used the notation
\begin{equation}
\pi_1(v)=(0,v), \quad
\pi_2(v,v')=v.
\end{equation}
As above we conclude that the integral operator
\begin{equation}
(\7G f)(x):=\int G(x,y) f(y)\, dy
\end{equation}
gives the inverse of $\7D^\dag\7D$.

The construction of the Green's function $F(x,y)$
and the corresponding integral operator $\7F$ giving the
inverse of $\6D^\dag\6D$ is completely analogous.

\subsection{Some formal calculus}
The calculations here are purely formal (symbolic) and apply to any symbols
with associative law of multiplication.
Let $F,F^{-1},G^{-1},Q_1,Q_2$ be such symbols.
We shall assume the relation
\begin{equation}\Label{}
F^{-1}=G^{-1}+Q_1Q_2.
\end{equation}

\bl\Label{formal1}
Suppose 
$$G^{-1}G=FF^{-1}=1.$$
Then $$(1-Q_2FQ_1)(1+Q_2GQ_1)=1, \quad (1-Q_2FQ_1)Q_2G=Q_2F.$$
\el

\bpf
We have
\begin{multline}
(1-Q_2FQ_1)(1+Q_2GQ_1)=1+Q_2GQ_1-Q_2FQ_1(1+Q_2GQ_1)\\
=1+Q_2GQ_1-Q_2F(1+Q_1Q_2G)Q_1
=1+Q_2GQ_1-Q_2F(G^{-1}G+Q_1Q_2G)Q_1\\
=1+Q_2GQ_1-Q_2F(G^{-1}+Q_1Q_2)GQ_1
=1+Q_2GQ_1-Q_2FF^{-1}GQ_1\\
=1+Q_2GQ_1-Q_2GQ_1=1
\end{multline}
for the first and
\begin{multline}\Label{}
(1-Q_2FQ_1)Q_2G=Q_2(1-FQ_1Q_2)G=Q_2(FF^{-1}-FQ_1Q_2)G\\
=Q_2F(F^{-1}-Q_1Q_2)G=Q_2FG^{-1}G=Q_2F
\end{multline}
for the second identity.
\epf

Analogously, we have the symmetric statement:

\bl\Label{formal2}
Suppose 
$$GG^{-1}=F^{-1}F=1.$$
Then $$(1+Q_2GQ_1)(1-Q_2FQ_1)=1, \quad GQ_1(1-Q_2FQ_1)=FQ_1.$$
\el

\bpf
The first identity can either be proved directly as above or
by exchanging $F,F^{-1}$ with $G,G^{-1}$ and changing the sign of $Q_1$
in the first identity of Lemma~\ref{formal1}.
For the second identity we have
\begin{multline}\Label{}
GQ_1(1-Q_2FQ_1)=G(1-Q_1Q_2F)Q_1=G(F^{-1}F-Q_1Q_2F)Q_1\\
=G(F^{-1}-Q_1Q_2)FQ_2=GG^{-1}FQ_2=FQ_2.
\end{multline}
\epf

\section{Instanton connection}\label{Solution}
Consider the kernel 
\begin{equation}\Label{kern}
\ker \6D^\dag\subset \G'_\L(S\otimes E) \oplus W_\L 
\end{equation}
and any linear isometry 
$$\Psi\colon {\6W}\to\ker\6D^\dag,$$
where $\6W$ is a hermitian vector space of the same dimension
as $\ker\6D^\dag$. 
We see $\Psi$ as an isometric parametrization of $\ker \6D^\dag$.
Then the adjoint operator $\Psi^\dag\colon \G'_\L(S\otimes E) \oplus W_\L\to \6W$
is the orthogonal projection to $\ker \6D^\dag$, composed with $\Psi^{-1}.$
In particular, we have
\begin{equation}\Label{ddag0}
\Psi^\dag \Psi = \id_{{\6W}}.
\end{equation}
Using the direct sum decomposition in \eqref{kern} we write
$$\Psi(w)=(\psi(w),\chi(w))\in \G'_\L(S\otimes E) \oplus W_\L$$
 for $w\in \6W$ and the condition $\Psi(w)\in \ker \6D^\dag$ means
 \begin{equation}\Label{ker-psi}
\7D^\dag\psi(w)+\7Q\chi(w)=0. 
\end{equation}
Assuming that $\7D^\dag$ is invertible, \eqref{ker-psi} can be rewritten as
\begin{equation}\Label{psichi}
\psi(w) = - (\7D^\dag)^{-1}\7Q\chi(w).
\end{equation}
Thus $\Psi$ is completely determined by its component $\chi$
and the latter defines a linear isomorphism $\chi\colon \6W\to W_\L$.
Furthermore, we have
\begin{equation}\Label{psidag}
\Psi^\dag(f,v) = \psi^\dag f+ \chi^\dag v, \quad (f,v)\in \G'_\L(S\otimes E) \oplus W_\L,
\end{equation}
and hence \eqref{ddag0} can be rewritten as
\begin{equation}
\chi^\dag \7Q^\dag \7D^{-1} (\7D^\dag)^{-1} \7Q\chi +\chi^\dag\chi = \id_{\6W}
\end{equation}
or, equivalently,
\begin{equation}\Label{chi-inv}
\7Q^\dag \7D^{-1} (\7D^\dag)^{-1} \7Q + \id_{V_\L} = (\chi^\dag)^{-1}\chi^{-1}.
\end{equation}

Consider $t$ varying in the set where $\7D$ and $\7D^\dag$ are invertible.
Over this set, the inverses of these operators depend smoothly on $t$.
Hence, $\ker\6D^\dag$ forms a subbundle of the trivial (infinite-dimensional) vector bundle 
\begin{equation}\Label{trivial}
\R^4\times (\G'_\L(S\otimes E) \oplus W_\L)\to \R^4.
\end{equation}
It has the natural {\em induced connection} $\2\nabla$ defined as the orthogonal projection
of the trivial connection on \eqref{trivial} to $\ker\6D^\dag$.

\bt[Nahm \cite{NahmCalorons}] The induced connection $\2\nabla$ is an {\em instanton connection}, i.e. the curvature of $\2\nabla$ is self-dual.
\et

Our goal here is to compute $\2\nabla$. Since $\ker\6D^\dag$ is a subbundle,
it has local trivializations over open subsets $U\subset\R^4$ given by linear maps 
$$\Psi=\Psi_t=(\psi_t,\chi_t)\colon W\to  \ker\6D_t^\dag, \quad t\in U,$$
where $\Psi_t$ is a linear isometry depending smoothly on $t$.
In this trivialization, we have
\begin{equation}
\2\nabla=d + A
\end{equation}
or, in coordinates $t=(t_\mu)$,
\begin{equation}
\2\nabla_{\frac{\d}{\d t_\mu}}  = \frac{\d}{\d t_\mu} + A_\mu, \quad \mu=0,1,2,3.
\end{equation}
By the construction of $\2\nabla$,
\begin{equation}\Label{A}
A_\mu = \Psi^\dag \d_\mu \Psi.
\end{equation}
Differentiating \eqref{ddag0}, we obtain
\begin{equation}
(\d_\mu\Psi^\dag) \Psi + \Psi^\dag \d_\mu \Psi=0
\end{equation}
and hence \eqref{A} can be rewritten as
\begin{equation}\Label{A'}
2A_\mu = \Psi^\dag \d_\mu \Psi - (\d_\mu\Psi^\dag) \Psi.
\end{equation}
Using \eqref{psidag} and \eqref{psichi} we rewrite \eqref{A'} as
\begin{multline}\Label{A''}
2A_\mu = \psi^\dag(\d_\mu\psi)- (\d_\mu\psi^\dag)\psi + \chi^\dag(\d_\mu\chi) -(\d_\mu \chi^\dag)\chi\\
=\chi^\dag\7Q^\dag \7D^{-1} \d_\mu((\7D^\dag)^{-1}\7Q\chi)
-\d_\mu(\chi^\dag\7Q^\dag \7D^{-1} )(\7D^\dag)^{-1}\7Q\chi
+ \chi^\dag(\d_\mu\chi) -(\d_\mu \chi^\dag)\chi\\
=\chi^\dag\7Q^\dag \big(\7D^{-1} \d_\mu (\7D^\dag)^{-1}
- (\d_\mu\7D^{-1}) (\7D^\dag)^{-1}\big)\7Q \chi\\
+ \chi^\dag \big(\id_{V_\L} + \7Q^\dag\7D^{-1}(\7D^\dag)^{-1}\7Q\big)(\d_\mu\chi)
- (\d_\mu\chi^\dag) \big(\id_{V_\L} + \7Q^\dag\7D^{-1}(\7D^\dag)^{-1}\7Q\big)\chi\\
=\chi^\dag\7Q^\dag \big(\7D^{-1} \d_\mu (\7D^\dag)^{-1}
- (\d_\mu\7D^{-1}) (\7D^\dag)^{-1}\big)\7Q \chi
+\chi^{-1}(\d_\mu\chi)-(\d_\mu\chi^\dag)(\chi^\dag)^{-1}
\end{multline}
where we have used \eqref{chi-inv} and the fact that $\7Q$ is independent of $t$.

We now compute the expression in the brackets in the first term of the last line in \eqref{A''}:
\begin{multline}\Label{brack}
\7D^{-1} \d_\mu (\7D^\dag)^{-1}
- (\d_\mu\7D^{-1}) (\7D^\dag)^{-1}
=\\
=-\7D^{-1}(\7D^\dag)^{-1}(\d_\mu\7D^\dag)(\7D^\dag)^{-1}+
\7D^{-1}(\d_\mu \7D)\7D^{-1}(\7D^\dag)^{-1}\\
= (\7D^\dag\7D)^{-1}\big( \7D^\dag(\d_\mu \7D) - (\d_\mu\7D^\dag)\7D\big) (\7D^\dag\7D)^{-1},
\end{multline}
where we have used the formula $\d_\mu(\7D^{-1})=-\7D^{-1}(\d_\mu\7D)\7D^{-1}$
and the analogous formula for $\7D^\dag$.
Setting 
$$G=\7G:=(\7D^\dag\7D)^{-1}, \quad F=\7F:=(\6D^\dag\6D)^{-1}, \quad Q_1:=\7Q,
\quad Q_2:=\7Q^\dag,$$
in Lemmas~\ref{formal1} and \ref{formal2}, we obtain
\begin{equation}\Label{formrel}
(\id-\7Q^\dag\7F\7Q)\7Q^\dag (\7D^\dag\7D)^{-1} = \7Q^\dag\7F, \quad
(\7D^\dag\7D)^{-1}\7Q (\id-\7Q^\dag\7F\7Q) = \7F\7Q.
\end{equation}
The operator $\7Q^\dag\7F\7Q\colon W_\L\to W_\L$ is a component of 
$$\6D\7F\6D^\dag\colon \G'_\L\oplus W_\L\to \G'_\L\oplus W_\L.$$
The latter operator satisfies $(\6D\7F\6D^\dag)^2=\6D\7F\6D^\dag$
and $\ker(\6D\7F\6D^\dag)=\ker\6D^\dag$
since $\7F$ is invertible and 
{$\ker\6D=0$ (see Section~\ref{ker}).}
Consequently, $\6D\7F\6D^\dag$ is the orthogonal projection
onto the orthogonal complement to $\ker\6D^\dag$
and $\id-\6D\7F\6D^\dag$ is the orthogonal projection onto $\ker\6D^\dag$.
On the other hand, $\Psi\Psi^\dag$ is also the orthogonal projection onto $\ker\6D^\dag$.
By uniqueness, we have 
\begin{equation}\Psi\Psi^\dag=\id-\6D\7F\6D^\dag\colon  \G'_\L\oplus W_\L\to \G'_\L\oplus W_\L.
\end{equation}
Taking suitable components of both sides, we obtain
\begin{equation}\label{chichi}
\chi\chi^\dag=\id-\7Q^\dag\7F\7Q,
\end{equation}
where $\7F=1\otimes F.$ 
In particular, $\id-\7Q^\dag\7F\7Q$ is invertible and 
we can rewrite \eqref{formrel} as
\begin{equation}\Label{rel'}
\7Q^\dag (\7D^\dag\7D)^{-1} = (\chi\chi^\dag)^{-1}\7Q^\dag\7F, \quad
(\7D^\dag\7D)^{-1}\7Q = \7F\7Q  (\chi\chi^\dag)^{-1}.
\end{equation}
Now using \eqref{brack} and \eqref{rel'}, we obtain
\begin{multline}\Label{brack'}
\7Q^\dag\big(\7D^{-1} \d_\mu (\7D^\dag)^{-1}
- (\d_\mu\7D^{-1}) (\7D^\dag)^{-1}\big)\7Q\\
=   (\chi\chi^\dag)^{-1}\7Q^\dag\7F\big( \7D^\dag(\d_\mu \7D) - (\d_\mu\7D^\dag)\7D\big)
\7F\7Q  (\chi\chi^\dag)^{-1} .
\end{multline}

We now use the explicit formulae \eqref{D} and \eqref{D-adj} to calculate the derivatives:
\begin{equation}
\d_\mu \7D= e_\mu\otimes\id, \quad \d_\mu \7D^\dag= \bar{e}_\mu\otimes\id,
\end{equation}
where $\bar{e}_0=e_0=\id$, $\bar e_j = e^\dagger_j= -e_j$ is the quaternion conjugation. 
We also consider decompositions similar to \eqref{qq}:
\begin{equation}
\7D=\sum_\nu e_\nu\otimes\7D_\nu, \quad 
\7D^\dag=\sum_\nu \bar{e}_\nu\otimes\7D_\nu,
\end{equation}
which are made explicit by \eqref{D} and \eqref{D-adj},
i.e.\ $$\7D_0=i\nabla+t_0, \quad \7D_j=T_j+t_j.$$
Then
\begin{equation}\Label{comm}
\7D^\dag(\d_\mu \7D) - (\d_\mu\7D^\dag)\7D
=\sum_\nu\bar e_{[\nu}e_{\mu]} \otimes \7D_\nu,
\end{equation}
where we use the notation 
$\bar e_{[\mu}e_{\nu]}:=\bar e_\mu e_\nu -\bar e_\nu e_\mu$.
On the other hand, we can use \eqref{ddqq'} to obtain
\begin{equation}\Label{dmu}
2\7D_\nu = \d_\nu (\6D^\dag\6D)=\d_\nu F^{-1},
\end{equation}
where $F\colon \G''_\L(E)\to\G_\L(E)$ is determined by $\7F=\id\otimes F$.
Subsituting \eqref{dmu} into \eqref{comm} and subsequently into \eqref{brack'}, we obtain
\begin{multline}\Label{brack''}
\7Q^\dag\big(\7D^{-1} \d_\mu (\7D^\dag)^{-1}
- (\d_\mu\7D^{-1}) (\7D^\dag)^{-1}\big)\7Q\\
=   \frac12\sum_\nu(\chi\chi^\dag)^{-1}\7Q^\dag (\id\otimes F)
\big( \bar e_{[\nu} e_{\mu]}\otimes\d_\nu F^{-1}  \big)
(\id\otimes F)\7Q  (\chi\chi^\dag)^{-1}\\
=-\frac12\sum_\nu(\chi\chi^\dag)^{-1}
\d_\nu\big(\7Q^\dag 
( \bar e_{[\nu} e_{\mu]}\otimes F )
\7Q\big)  (\chi\chi^\dag)^{-1},
\end{multline}
where we have used the fact that $\7Q$ is independent of $t$
and the formula $\d_{\mu}F^{-1}=-F^{-1}(\d_{\mu}F) F^{-1}$.
Using the Green's function $F(x,y)$ of the operator $F$ and the explicit formulae \eqref{Qdag} and \eqref{Q}, we calculate
\begin{multline}\Label{qfq}
\big(\7Q^\dag 
( \bar e_{[\nu} e_{\mu]}\otimes F )
\7Q\big)(v_1,\ldots,v_n)=\\
=\Big(Q_\beta^\dag \int \big(\bar e_{[\nu} e_{\mu]}\otimes F(\l_\beta,s)\big) \sum_\alpha \delta(s-\l_\alpha) Q_\alpha v_\alpha \, ds\Big)_{1\le \beta\le n}\\
=\Big(\sum_\alpha Q_\beta^\dag \big(\bar e_{[\nu} e_{\mu]}\otimes F(\l_\beta,\l_\alpha)\big)   Q_\alpha v_\alpha\Big)_{1\le \beta\le n}=\\
=\Big(Q^\dag \big(\bar e_{[\mu} e_{\nu]}\otimes F({\l_*},\l_*)\big) Q\Big) (v_1,\ldots,v_n),
\end{multline}
where we have used the notation
$$Q\colon W_\L\to (S\otimes E)_\L:=\oplus_{\alpha=1}^n (S\otimes E)_{\l_\alpha},\quad
(v_\alpha)_{1\le \alpha\le n}\mapsto (Q_\alpha v_\alpha)_{1\le \alpha\le n},$$
and 
$$F({\l_*},\l_*)\colon E_\L\to  E_\L:=\oplus_{\alpha=1}^n E_{\l_\alpha}, \quad
(w_\alpha)_{1\le \alpha\le n}\mapsto \Big(\sum_\alpha F(\l_\beta,\l_\alpha)w_\alpha\Big)_{1\le \beta\le n}.$$
Substituting \eqref{qfq} into \eqref{brack''} and subsequently into \eqref{A''}
we finally obtain
\begin{multline}\Label{A'''}
2A_\mu 
=-\frac12\sum_\nu
\chi^{-1} Q^\dag \big(\bar e_{[\nu} e_{\mu]}\otimes \d_\nu F({\l_*},\l_*)\big) Q (\chi^\dag)^{-1}\\
+\chi^{-1}(\d_\mu\chi)-(\d_\mu\chi^\dag)(\chi^\dag)^{-1}.
\end{multline}

To summarize, one can use  \eqref{chichi} to find $\chi,$ for example by choosing $\chi=\chi^\dag=\left(\id-\7Q^\dag\7F\7Q\right)^{\frac{1}{2}},$ and the above formula \eqref{A'''} to obtain the caloron connection.

\section{Conclusions}
Current resurgence of interest in instantons on various spaces is due to their significance in geometry, gauge theory, string theory,  and representation theory.  There are three natural directions in which our result can be generalized.  One would like to allow for arbitrary monopole charges for the caloron.  This would amount to working with Nahm data that changes rank across the $\lambda$-points, thus having natural projection operators at such points.  It should be a fairly straightforward exercise incorporating such projectors into our formula.  The second generalization involves considering instantons on curved spaces, such as a hyperk\"ahler ALF space, for example.  The Nahm data in this case is given in terms of a bow \cite{Cherkis:2010bn} and our formula allows for a generalization to this case as well \cite{Clare}.  Lastly, one would like to have an expression for a caloron connection for arbitrary gauge group $G$.  This remains an open problem.

\end{document}